\documentclass[12pt,titlepage]{article}
\usepackage{amsmath}
\usepackage{amssymb}
\usepackage{epic,eepic}
\usepackage{verbatim}

\def\R{\mathbb{R}}

\def\al{\alpha}
\def\bt{\beta}
\def\lm{\lambda}
\def\norm#1{\lVert#1\rVert}
\def\abs#1{\lvert#1\rvert}
 
\def\eps{\varepsilon}

\def\TGS{T_{GS}}
\def\TGSJ{T_{GSJ}}
\def\VGS{V_{GS}}
\def\VGSJ{V_{GSJ}}

\def\dV{\partial{V}}

\usepackage{amsthm}
\newtheorem{lem}{Lemma}[section]
\newtheorem{thm}{Theorem}[section]

\newtheorem{exl}{Example}

\theoremstyle{remark}
\newtheorem{rem}{Remark}[section]

\renewcommand{\baselinestretch}{2.0}

\begin{document}
\title{Acceleration Operators in the Value Iteration Algorithms for Markov Decision
Processes \\} \vskip5mm
\author{Oleksandr Shlakhter$^1$, Chi-Guhn Lee$^{1,}$\footnotemark[2],
Dmitry Khmelev\footnotemark[3], Nasser Jaber$^1$ \\ \\
$^1$Department of Mechanical and Industrial Engineering \\
$^3$Department of Mathematics \\
University of Toronto, Toronto, Ontario, Canada,  M5S 3G8}
\footnotetext[2]{Corresponding author can be reached at
\texttt{cglee@mie.utoronto.ca}} \footnotetext[3]{Passed away in 2004 while our
research was in progress} \maketitle

\begin{abstract}
We study the general approach to accelerating the convergence of the most widely used
solution method of Markov decision processes with the total expected discounted
reward. Inspired by the monotone behavior of the contraction mappings in the feasible
set of the linear programming problem equivalent to the MDP, we establish a class of
operators that can be used in combination with a contraction mapping operator in the
standard value iteration algorithm and its variants. We then propose two such
operators, which can be easily implemented as part of the value iteration algorithm
and its variants. Numerical studies show that the computational savings can be
significant especially when the discount factor approaches 1 and the transition
probability matrix becomes dense, in which the standard value iteration algorithm and
its variants suffer from slow convergence.

\end{abstract}

\section{Introduction} \label{sec:intro}
Consider an infinite horizon Markov decision process (MDP) with a finite set of $i$
states denoted by $S$, a finite set of actions $A(i)$ for each state $i \in S$, an
immediate reward $r(i,a)$ for each $i \in S$ and $a \in A=\cup_{i \in S}A(i)$, and a
transition probability $p_{ij}(a)$ for each $i,j \in S$ and $a\in A(i)$. The
objective is to determine $v_i$, the maximum expected total discounted reward over an
infinite horizon starting in state $i$, where $\lambda$ is the discount factor ($0\le
\lambda < 1$). It is well known~\cite{Puterman:MDP94} that $v$ satisfies the
optimality equation
\begin{equation} \label{MDP}
v_i=\max_{a\in A(i)}\left\{r(i,a)+\lambda \sum_{j\in S} p_{ij}(a) v_j\right\},
\end{equation}
and the actions attaining the maximum in Equation (\ref{MDP}) give rise to an optimal
stationary policy. Let $U$ denote the set of bounded real valued functions on $S$
with the norm $||v||=\max_{i\in S}|v_i|$, and let the function $d:S\rightarrow A(s)$ specifies
the action choice for state $i\in S$ (i.e., $d(i) \in A(i)$ for each $i \in S$).
Furthermore, let $r_d$ denote the $|S|$-vector, with $i$-th component $r_d(i)\equiv
r(i,d(i))$, and $P_d$ the $|S|\times|S|$ matrix with $(i,j)$-th entry given by
$p_d(j|i)\equiv p_{ij}(d(i))$. Then the optimality equation given in Equation
(\ref{MDP}) can be written, with the definition of the operator $T$ on $U$, in the
following vector notation:
\begin{equation} \label{MDP vector}
v=Tv\equiv \max_{d\in \Pi}\left\{r_d + \lambda P_d v \right\},
\end{equation}
where $\Pi$ is the set of policies.

There are several standard methods for finding optimal or approximately optimal
policies for the MDP. Approaches widely employed to solve MDP problems include value
iteration and policy iteration~\cite{Puterman:MDP94}. Although simple to implement,
these approaches are nevertheless limited in the size of problems that can be solved
due to the excessive computation required to find close-to-optimal solutions.
Techniques to improve the convergence of the value iteration algorithm have been
studied in~\cite{Herzberg/Yechiali:EJOR96p622,Herzberg/Yechiali:OR94}, just to name a
few. The standard value iteration algorithm is summarized as
follows~\cite{Puterman:MDP94}

\renewcommand{\baselinestretch}{1.5}
\noindent{\bf Value Iteration (VI) Algorithm} \vspace{-0.3in}
\begin{quotation}
\texttt{
\begin{description}
  \item [Step 0] Select $v^0\in\R^I$, set $n=0$, and specify $\eps >0$.
  \item [Step 1] Compute $v^{n+1}_i=Tv^{n}_i$  for all $i\in S$.
  \item [Step 2] If $\norm{v^{n+1}-v^{n}}<\eps(1-\lm)/2\lm$, go to {\bf Step 3}.
  Otherwise, increase $n$ by 1 and return to {\bf Step 1}.
  \item [Step 3] Return with the actions attaining the maximum in {\bf Step 1}.
\end{description}
}
\end{quotation}
\renewcommand{\baselinestretch}{2.0}

An alternative approach is based on linear programming. Although linear
programming-based approaches are generally dismissed as
inefficient~\cite{Puterman:MDP94}, they have recently inspired the research community
to find new approaches using the well-studied theory of linear programming. Widely
employed linear programming
formulations~\cite{DeFarias/VanRoy:OR03p850,Puterman:MDP94,Trick/Zin:MD97p255} are
due to Derman~\cite{Derman:MDP70} as given below:
\begin{equation} \label{eq:LP}
\min \left\{\sum_i v_i \mid{} v_i-\lm\sum_jp_{ij}(a)v_j\ge{}r(i,a), \forall i\in{S},
\forall a\in A(i), v \in \R^{|S|} \right\}.
\end{equation}
Notice that the LP is a minimization problem as opposed to the given MDP being a
maximization problem. This is due to the implementation of the contraction mapping as
constraints of the LP.

 One of the recent methods rooted in the above linear programming
formulation is an approximation approach using basis
functions~\cite{DeFarias/VanRoy:OR03p850}. They use a linear combination of
pre-selected basis functions to fit the value function of MDP to obtain an approximate
value function, which is then optimized in the feasible set of the linear program --
named the approximate linear program. While the quality of the approximation
heavily depends on the pre-selected basis functions, they fail to provide reliable
guidelines on how to chose basis functions. The feasible set $V$ of the equivalent
linear program can be written in vector notation as follows:
\begin{equation}
  V=\{v\in\R^{|S|} \mid v\ge{}Tv\}. \label{eq:VJ}
\end{equation}

We propose a new class of operators that can be integrated into the standard and
variants of the value iteration algorithm for Markov decision processes so as to
speed up the convergence of the iterative search. The speed of convergence will be
measured in both the number of iterations and the CPU time. The characterization of
the operator class is motivated by the monotone behavior of contraction mappings in
set $V$ and the ideas of the basis function approximation. The major contribution of
this paper is the identification of a new class of operators which accelerate the
value iteration algorithm and its variants. The observation of monotone behavior of
the contraction mapping in set $V$ in Equation (\ref{eq:VJ}) would also inspire the
research community to explore other types of operators.

The rest of the paper is organized as follows. Section~\ref{sec:accel-value}
introduces conditions which need to be satisfied by accelerating operators, that are
to be combined with value iteration operators. Section~\ref{sec:numerical-studies}
presents the performance improvements due to the acceleration as well as a discussion
on computational complexity of the acceleration operators. In section
\ref{sec:total-reward} we discuss the application of the acceleration operators to
the expected total reward, semi-Markov and continuous-time Markov decision processes.
Finally, section~\ref{sec:conclusions} concludes the paper.

\section{Accelerated Value Iterations Algorithms}
\label{sec:accel-value}

In this section we will introduce a class of operators that will speed up the
convergence of the standard value iteration and its variants. Operators in this class
can be used in combination with the contraction mapping of the value iteration
algorithm and its variants. We will also present two such operators later in this
section.

The most crucial observation, that leads to characterization of the acceleration
operators, is that the operator $T$ is a component-wise monotone mapping in $V$ and
the set $V$ is invariant under $T$ as formally stated in Lemma \ref{lem:invariance}.

\begin{lem}
  \label{lem:invariance}
  $V$ is invariant under $T$. That is, $TV \subset V$.
\end{lem}

It is also known~\cite{Puterman:MDP94} that the optimal solution $v^*$ of the MDP or
the fixed point of the operator $T$ given in (\ref{MDP vector}) is such that
\[
v^* \le v, \; \forall v \in V.
\]

Lemma~\ref{lem:invariance} suggests the following conditions for acceleration
operator $Z$ on the set $V$.

\bigskip

\noindent {\bf Acceleration Conditions}
\begin{description}
\item [(A)] $Z V\subset V$,
\item [(B)] $Z v\le v, \; \forall v \in V$.
\end{description}
\bigskip

One can see that the identity operator and the operator $T$ satisfy both of the
conditions. For a given operator satisfying two conditions {\bf (A)} and {\bf (B)},
an accelerated value iteration algorithm can be defined. Below is a summary of the
Generalized Accelerated Value Iteration Algorithm (GAVI) proposed in the present
paper.

\bigskip
\renewcommand{\baselinestretch}{1.5}
\noindent{\bf General Accelerated Value Iteration Algorithm (GAVI)} \vspace{-0.4in}
\begin{quotation}
\texttt{
\begin{description}
  \item [Step 0] Select $w^0\in{}V$, set $n=0$, and specify $\eps >0$.
  \item [Step 1] Compute $w^{n+1}_i=Z Tw^{n}_i$ for all $i\in{S}$.
  \item [Step 2] If $\norm{Tw^{n}-w^{n}}\le \eps(1-\lm)/2\lm$, go to {\bf Step 3}.
  Otherwise increase $n$ by 1 and return to {\bf Step 1}.
  \item [Step 3] Return with the actions attaining the maximum in {\bf Step 1}.
\end{description}
}
\end{quotation}
\renewcommand{\baselinestretch}{2.0}

Notice that different acceleration operators may be used in different iterations of
the value iteration algorithm, in which case $Z_n$ is an acceleration operator used
in Step 1 of $n$-th iteration, instead of $Z$. In this paper we restrict our
attention to the case where $Z_n=Z$ for all $n$.

\begin{rem} \label{r:in.-point}
GAVI must be initialized with $w^0\in V$. It is shown
in~\cite{Schweitzer/Seidmann:JMAA85p568} that there always exists a feasible solution
$u (\in V)$ of the equivalent LP in form $u(i)=\alpha \vec{1}$ for all $i$, where
$\vec{1}$ is a vector of ones. Indeed, taking $\alpha=(1/(1-\lambda))\max_{i\in
S,a\in A(i)}r(i,a)$, we will have $u$ be a feasible solution of the above LP.
\end{rem}

\begin{rem}
 \label{r:convergence}
It is easy to show that if $v^*$ is the fixed point of operator $T$, $v^{n+1}=
Tv^{n}$, and $w^0=v^0$. Then
  \begin{enumerate} \vspace{-0.2in}
  \renewcommand{\labelenumi}{(\roman{enumi})}
  \item $v^*\le{}w^n\le{}v^n$,\vspace{-0.2in}
  \item $v^*=\lim_{n\to\infty}w^n=\lim_{n\to\infty}v^n$,\vspace{-0.2in}
  \item $\norm{v^*-w^n}\le\norm{v^*-v^n}$.\vspace{-0.2in}
  \item The sequence $\{w^n\}$ converges globally on $V$(for all $w^0 \in V$) with order 1 at a rate less than or equal to $\lm$; its global asymptotic average rate of convergence is less than or equal to $\lm$$^1$.
  \end{enumerate}
\end{rem}

\footnotetext[1]{An extensive discussion of convergence orders and rates can be found in \cite{Puterman:MDP94}.}

As mentioned earlier, $T$ itself is such an operator, resulting in an algorithm that
performs {\bf Step 1} of VI twice per iteration. Trivially, such modification
accelerates the convergence of the algorithm but each iteration is twice as expensive
as VI. Therefore, it is important to identify an acceleration operator satisfying two
conditions {\bf (A)} and {\bf (B)} that requires little additional computation, so
that the reduction in the number of iterations before convergence is more dramatic than
the increased computation per iteration. If the condition {\bf (B)} holds strictly, the
improvement per iteration is guaranteed to be more than that of the standard value
iteration (VI).

Now we propose an acceleration operator that requires little additional computation
per iteration but reduces the number of iterations significantly. Understanding that
the value iteration algorithm is simply finding a component-wise minimum point in the
set $V$, any operator that maps a given point (vector) $v$ in $V$ to another vector
$u$ in $V$ such that $u\le v$ with the strict inequality with at least one component
will accelerate the convergence.

\bigskip

\noindent {\bf Projective Operator}

\noindent For $v \in V$, $Z:v \rightarrow \al^* v$, where $\al^*$ is the optimal
solution of the following optimization problem:
\begin{equation} \label{eq:1DLP}
\min \left\{\alpha \sum v_i \mid T (\alpha v) \le \alpha v \right\}.
\end{equation}

\bigskip

\begin{figure} [p]
  \begin{center}
    \setlength{\unitlength}{0.00083333in}
    {\renewcommand{\dashlinestretch}{30}
      \begin{picture}(3841,3700)(0,500)
        \thicklines
        \put(1882,1950){\blacken\ellipse{70}{70}}
        \put(1882,1950){\ellipse{70}{70}}
        \put(2967,2325){\blacken\ellipse{70}{70}}
        \put(2967,2325){\ellipse{70}{70}}
        \put(2217,2250){\blacken\ellipse{70}{70}}
        \put(2217,2250){\ellipse{70}{70}}
        \put(1429,1905){\blacken\ellipse{70}{70}}
        \put(1429,1905){\ellipse{70}{70}}
        \put(1332,1785){\blacken\ellipse{70}{70}}
        \put(1332,1785){\ellipse{70}{70}}
        \put(357,1635){\blacken\ellipse{70}{70}}
        \put(357,1635){\ellipse{70}{70}}
        \thinlines
        \blacken\path(72.000,3480.000)(42.000,3600.000)(12.000,3480.000)(72.000,3480.000)
        \path(42,3600)(42,225)(3742,225)
        \blacken\path(3622.000,195.000)(3742.000,225.000)(3622.000,255.000)(3622.000,195.000)
        \thicklines
        \path(1258,3216)(358,1641)(1408,1791)(3808,2616)
        \thinlines
        \dashline{60.000}(42,225)(3117,3075)
        \dashline{60.000}(34,225)(2337,2993)
        \dashline{60.000}(2217,2250)(2965,2333)
        \dashline{60.000}(1392,1883)(1880,1959)
        \put(-100,0){\makebox(0,0)[lb]{0}}
        \put(1362,2737){\makebox(0,0)[lb]{$V$}}
        \put(3042,2100){\makebox(0,0)[lb]{$w^{n-1}$}}
        \put(1955,2355){\makebox(0,0)[lb]{$Tw^{n-1}$}}
        \put(1849,1672){\makebox(0,0)[lb]{$w^n$}}
        \put(1167,2033){\makebox(0,0)[lb]{$Tw^n$}}
        \put(192,1410){\makebox(0,0)[lb]{$v*$}}
        \put(972,1822){\makebox(0,0)[lb]{$w^{n+1}$}}
        \put(-200,3500){\makebox(0,0)[lb]{$v_2$}}
        \put(3700,0){\makebox(0,0)[lb]{$v_1$}}
      \end{picture}
    }
  \end{center}
  \caption{Projective AVI}
  \label{f:PAVI}
\end{figure}
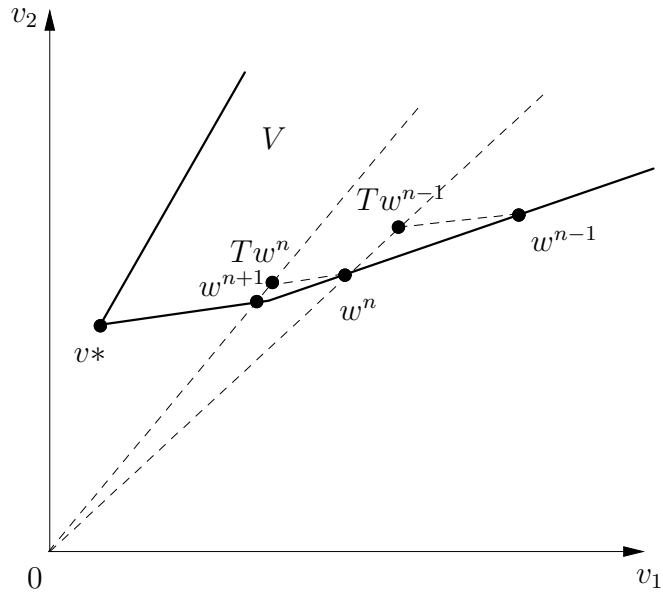

We will call GAVI with Projective Operator as Projective Accelerated Value Iteration
or PAVI for short in the paper. It is interesting to view Projective Operator as an
approximation approach using a basis as is proposed
in~\cite{DeFarias/VanRoy:OR03p850} since the optimization problem given in (\ref{eq:1DLP})
is nothing but the approximate LP with a single basis function $v$ found by the
standard operator $T$ in {\bf Step 1} of VI or GAVI. The choice of basis functions is
not an issue in this approach since $V$ is invariant under $T$ and $T$ is a
contraction mapping. That is, $T$ finds a vector $v \in V$ so that optimizing the
approximate linear program with a single decision variable, denoted by $\alpha$ in
the definition of Projective Operator, results in a new  vector closer to the fixed
point. The role of Projective Operator is graphically
illustrated in Figure~\ref{f:PAVI}, where $Z$ projects the given point $v \in V$ to
the boundary of $V$.

\vspace{1in}
\begin{center}
Figure~\ref{f:PAVI} goes here.
\end{center}
\vspace{1in}

\begin{thm} \label{l:Pu_satisfy12}
Projective Operator $Z$ satisfies the conditions {\bf (A)} and {\bf (B)} if
$r(i,a)\ge 0$ for all $i\in S$, for all $a\in A(i)$.
\end{thm}

\begin{rem}
 \label{r:non-negativity}
The non-negativity condition is not really restrictive. The rewards can be adjusted
to be non-negative for all states and actions
\[
   r(i,a)'=r(i,a)+{}\max_{j\in{}S, a\in{}A(j)}\bigl\{|r(j,a)| \bigr\}, \; \forall i\in
   S,\forall a\in A(i).
\]
The new MDP shares the optimal policy with the original so that it will not matter
which MDP to solve. For the expected total discounted MDP, the reward adjustment can
also be considered as the coordinate translation as follows:
\[
    v_i'=v_i+{}\max_{j\in{}S, a\in{}A(j)}\bigl\{|r(j,a)| / (1-\lm) \bigr\}, \;
    \forall i\in S,
\]

\end{rem}

Let $\dV$ denote the set of the boundary points of $V$ and $int(V)$ the set of the
interior points of $V$, i.e., $int(V)=V\setminus\dV$.

\begin{lem}
  \label{lem:fully dense}
Suppose that $P_d$ is fully dense (no zero entry) for all $d\in \Pi$, then $Tv < v$
for all $v\in V$ with $v \neq v^*$.
\end{lem}

\begin{rem}
Notice that if the transition probability matrix $P_d$ is fully dense for all $d\in
\Pi$, the strict inequality in (i) and (iii) of the Remark \ref{r:convergence} is guaranteed, i.e.
the sequence of iterates PAVI converges faster than that of VI. This is because it is always true
that for all $v\in V$ with $v \neq v^*$, $Tv \in int(V)$ and Projective Operator maps
$Tv \in int(V)$ to a vector $w \in \dV$. It is obvious that $w < v$. It worth to mention, however, that strict inequalities in (i) and (iii) do not guaranty that the stopping criteria of value iteration algorithm will be satisfied strictly faster.
\end{rem}

\begin{lem} \label{lem:interior}
The interior of $V$ is invariant under $T$. That is,
\[
T(int(V)) \subset int(V).
\]
\end{lem}

Lemma~\ref{lem:interior} suggests a variation of Projective Operator for which the strict inequality in (i) and (iii) of the Remark \ref{r:convergence} is guaranteed, which allows
the sequence of iterates PAVI to converge faster than that of VI in sense of (i) and (iii). Given a factor $\bt (0<\bt<1)$, one
can define a variation of Projective Operator as $Z^\bt{u}=(1-\bt)Zu+\bt{u}$, which
always finds an interior point between $Tw^n$ and $ZTw^n=w^{n+1}$ so that in the next
iteration $Tw^{n+1}$ is guaranteed to be an interior point of $V$.

Now we present a new acceleration operator that satisfies {\bf Acceleration
Conditions} {\bf (A)} and {\bf (B)}.

\bigskip

\noindent {\bf Linear Extension Operator}

\noindent For $v\in V$, $Z:v\rightarrow v+\al^*(u-v)$, where $u=Tv$ and $\al^*$ is
the optimal solution to the following optimization problem:
\begin{equation} \label{eq:1DLP2}
\min\left\{\sum v_i + \alpha\sum (u_i-v_i) \mid T(v+\alpha (u-v)) \le v+\alpha
(u-v)\right\}.
\end{equation}

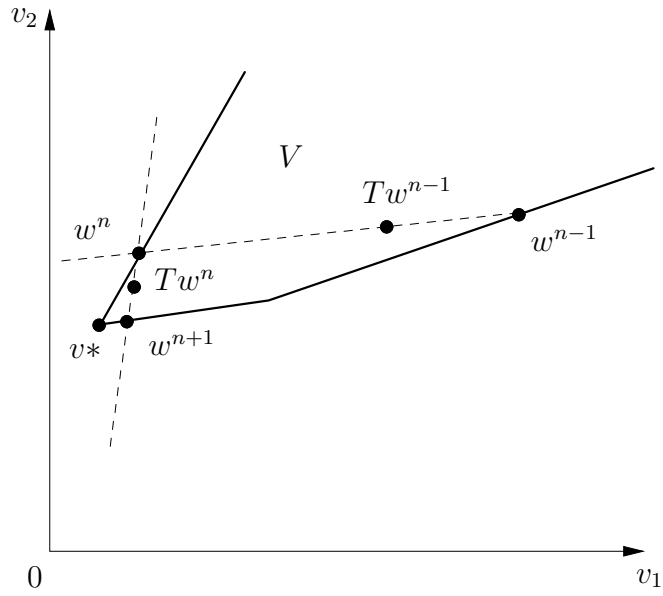
\begin{figure} [b]
\begin{center}
  \setlength{\unitlength}{0.00083333in}
  \renewcommand{\dashlinestretch}{30}
    \begin{picture}(3841,3264)(0,500)
      \thicklines
      \put(2967,2325){\blacken\ellipse{70}{70}}
      \put(2967,2325){\ellipse{70}{70}}
      \put(2142,2250){\blacken\ellipse{70}{70}}
      \put(2142,2250){\ellipse{70}{70}}
      \put(567,1875){\blacken\ellipse{70}{70}}
      \put(567,1875){\ellipse{70}{70}}
      \put(522,1658){\blacken\ellipse{70}{70}}
      \put(522,1658){\ellipse{70}{70}}
      \put(597,2085){\blacken\ellipse{70}{70}}
      \put(597,2085){\ellipse{70}{70}}
      \put(349,1635){\blacken\ellipse{70}{70}}
      \put(349,1635){\ellipse{70}{70}}
      \thinlines
      \blacken\path(72.000,3480.000)(42.000,3600.000)(12.000,3480.000)(72.000,3480.000)
      \path(42,3600)(42,225)(3742,225)
      \blacken\path(3622.000,195.000)(3742.000,225.000)(3622.000,255.000)(3622.000,195.000)
      \thicklines
      \path(1258,3216)(358,1641)(1408,1791)(3808,2616)
      \thinlines
      \dashline{60.000}(116,2037)(2965,2337)
      \dashline{60.000}(709,2933)(417,878)
      \put(-100,0){\makebox(0,0)[lb]{0}}
      \put(1467,2625){\makebox(0,0)[lb]{$V$}}
      \put(3042,2100){\makebox(0,0)[lb]{$w^{n-1}$}}
      \put(1992,2400){\makebox(0,0)[lb]{$Tw^{n-1}$}}
      \put(192,2175){\makebox(0,0)[lb]{$w^n$}}
      \put(694,1853){\makebox(0,0)[lb]{$Tw^n$}}
      \put(657,1448){\makebox(0,0)[lb]{$w^{n+1}$}}
      \put(162,1433){\makebox(0,0)[lb]{$v*$}}
      \put(-200,3500){\makebox(0,0)[lb]{$v_2$}}
      \put(3700,0){\makebox(0,0)[lb]{$v_1$}}
    \end{picture}
\end{center}
\caption{Linear AVI} \label{f:LAVI}
\end{figure}

Figure \ref{f:LAVI} graphically illustrates how Linear Extension Operator works. It
casts $Tw^n$ in the direction of $Tw^n-w^n$ to the boundary of set $V$. For $w^n \in
V$, we have $Tw^n \le w^n$, implying that $Tw^n-w^n$ is an improving direction within
set $V$. As a result, Linear Extension Operator moves $Tw^n$ closer to the fixed
point. An improvement is guaranteed whenever $Tw^n \in int(V)$. When the
transition matrix $P_d$ is fully dense for all $d \in \Pi$, such a improvement is also
guaranteed. We can also define a variation of Linear Extension Operation as we did
with Projective Operator to guarantee a strict inequalities in (i) and (iii) of the Remark \ref{r:convergence} even when the
matrices are not fully dense.

\vspace{1in}
\begin{center}
Figure~\ref{f:LAVI} goes here.
\end{center}
\vspace{1in}

\begin{thm}\label{thm:L_w_satisfy12}
Linear Extension Operator $Z$ satisfies the conditions {\bf (A)} and {\bf (B)}.
\end{thm}

When Linear Extension Operator  is used in {\bf Step 1} of GAVI, we call the
algorithm as Linear Extension Accelerated Value Iteration or LAVI for short in the
paper.

In what follows we will combine the two acceleration operators - Projective Operator
and Linear Extension Operator - with variants of the standard value iteration
algorithms: Jacobi, Gauss-Seidel, and Gauss-Seidel-Jacobi value iteration algorithms.

\noindent \textbf {Jacobi:} $v^{n+1}=T_Jv^n$ where
\begin{equation}
  v^{n+1}_i=\max_{a\in{A(i)}}\left\{\Bigl[ r(i,a)+\lm\sum_{j \neq i} p_{ij}(a)
  v^n_j \Bigr ]/ {\Bigl [ 1-\lm p_{ii}(a)\Bigr ]}\right\}, \; \forall i\in S. \label{eq:TJ}
\end{equation}
\textbf {Gauss-Seidel:} $v^{n+1}=\TGS{v^n}$ where
\begin{equation}
  v^{n+1}_i=\max_{a\in{A(i)}}\left\{r(i,a)+\lm\sum_{j<i} p_{ij}(a)
  v^{n+1}_j+ \lm\sum_{j \geq i} p_{ij}(a) v^n_j \right\}, \; \forall i\in S. \label{eq:TGS}
\end{equation}
\textbf {Gauss-Seidel-Jacobi:} $v^{n+1}=\TGSJ{v^n}$ where
\begin{equation}
  v^{n+1}_i=\max_{a\in{A(i)}}\left\{ \Bigl[ r(i,a)+\lm\sum_{j<i} p_{ij}(a)
  v^{n+1}_j+ \lm\sum_{j>i} p_{ij}(a) v^n_j \Bigr ]/
              {\Bigl [ 1-\lm p_{ii}(a)\Bigr ]}\right\}, \; \forall i\in S. \label{eq:TGSJ}
\end{equation}
Notice that the definition of the operators $\TGS$ and $\TGSJ$ depends on the order
of state index. We will show that these operators can be used in the position of the
standard operator $T$ in GAVI (and hence in PAVI and in LAVI). We start with the
following definition of sets:
\begin{align*}
V_J&=\{v\in\R^{|S|}\mid v\ge{}T_Jv\}, \\
\VGS&=\{v\in\R^{|S|}\mid v\ge{}\TGS{v}\}, \\
\VGSJ&=\{v\in\R^{|S|}\mid v\ge{}\TGSJ{v}\}.
\end{align*}

The following lemma is an analogue of Lemma~\ref{lem:invariance}:
\begin{lem}
  \label{lem:invariance of splittings}
$V_J$, $V_{GS}$, and $V_{GSJ}$ are invariant under $T_J,T_{GS}$, and $T_{GSJ}$
respectively.
\begin{comment}
  \noindent \begin{enumerate}
  \renewcommand{\labelenumi}{(\roman{enumi})}
  \item $V_J$ is invariant under $T_J$,
  \item $V_{GS}$ is invariant under $T_{GS}$,
  \item $V_{GSJ}$ is invariant under $T_{GSJ}$.
  \end{enumerate}
\end{comment}
\end{lem}

With Lemma \ref{lem:invariance of splittings} acceleration operators satisfying
conditions {\bf (A)} and {\bf (B)} can be used in {\bf Step 1} of GAVI with the
variants $T_J,T_{GS}$, and $T_{GSJ}$ of the standard operator $T$. However, it is not
trivial to define $V_{GS}$ and $V_{GSJ}$ with a set of linear inequalities and the
acceleration operators proposed in this paper will not work. To avoid the
problem, we should restrict the acceleration operators to a strict subset of $V_{GS}$
and $V_{GSJ}$. In this paper we use $V$ since it is a strict subset of $V_{GS}$ and
$V_{GSJ}$ as shown in Lemma \ref{l:relationships_between_V's} below.

\begin{lem} \label{l:relationships_between_V's}
The following relations hold
\begin{equation}
   V=V_J, \;\;\; \VGS=\VGSJ,\; \text{ and } \;  V\subset{\VGS}.\label{eq:3}
\end{equation}
\end{lem}

\begin{rem}
Gauss-Seidel methods require special consideration, since in general $V\ne{}\VGS$.
Indeed, let's consider an MDP with two states and one action per state with a
transition probability matrix $\genfrac(){0pt}1{0\,\, 1}{1\,\, 0}$, and a reward
vector $r=\genfrac(){0pt}1{1}{1}$.  Then $v \ge Tv$ is equivalent to $v_1\ge1+\lm
v_2$, $v_2\ge1+\lm v_1$, while $v \ge T_{GS}V$ is equivalent to $v_1\ge1+\lm v_2$,
$v_2\ge1+\lm(1+\lm v_1)$, and it is easy to see that these systems are different.
\end{rem}

\begin{lem} \label{l:TX(VX)subsetV}
\begin{equation}
    \TGS( \VGS) \subset V,\;\;\; \text{and} \;\;\;\; \TGSJ(\VGSJ)\subset V.
\end{equation}
\end{lem}

\begin{thm} \label{thm:V_invariance_under_TJ,TGS,TGSJ}
The set $V$ is invariant under $T_J$, $\TGS$, and $\TGSJ$. That is,
\[
  T_J{V}\subset{}V, \;\;\; \TGS{V}\subset{}V, \;\;\; \TGSJ{V}\subset{V}.
\]
\end{thm}

Theorem \ref{thm:V_invariance_under_TJ,TGS,TGSJ} states that $V$ is invariant under
$T_J,T_{GS}$, and $T_{GSJ}$, which suggests that these operators can replace $T$ in
GAVI to give rise to new accelerated value iteration algorithms. Therefore, we obtain
8 accelerated versions of the value iteration algorithm, which are conveniently
written in the form XAY, where `X' stands for either ``Projective'' or ``Linear
Extension'', `A' for ``Accelerated'', and `Y' for one of VI, J, GS, or GSJ. For example,
LAGS denotes Linear Extension `Accelerated Gauss-Seidel method with $w^{n+1}=ZT_{GS}$
for {\bf Step 1}. Non-accelerated versions will be shortened to VI, J, GS, and GSJ
without prefixes.

\begin{rem}
There are other modifications of VI based on the splittings of transition probability
matrix~\cite{Puterman:MDP94}. We believe that the strong inclusion demonstrated in
Lemma~\ref{l:TX(VX)subsetV} holds for the most practically useful splittings, but to
avoid unnecessary complications we do not prove the statement in the most general
setting.
\end{rem}

\section{Numerical Studies}
\label{sec:numerical-studies}

\subsection{Random examples}
\label{sec:random-example}

In this section we present numerical studies to demonstrate the computational
improvement that the two accelerating operators achieve over the standard value
iteration algorithm as well as its variants. We will consider two families of
randomly generated MDP problems. In all cases the number of actions in each state was
generated using a uniform random number generator between 2 and 99 and the immediate
reward for each state and action pair was generated using a uniform random number
generator in the range $(1, 100)$. Transition probability matrices were generated
differently between the two families. In Example 1, the transition probability matrices
have non-zero elements distributed uniformly, while in Example 2, non-zero elements
are only in a band around the diagonal.

In Example 1, we first fixed the number of non-zero entries in each row so that the
density of non-zero entries in that row should be equal to a given density level. We
randomly generated non-zero entries according to a uniform distribution over (0,1),
normalized these non-zero entries so that they add up to give the sum of 1, and then placed them
randomly across the row. In Example 2, non-zero entries were generated and placed
only around the diagonal so that the non-zero entries form a band around the
diagonal in the matrix. The number of non-zero entries in each row was determined for
a given density level.

We also used exponential distributions instead of uniform distributions but the numerical
results did not show any noticeable differences in performance. Therefore, we only
present numerical studies with uniform distributions in this paper.

In both Example 1 and 2, the density of the non-zero elements in the matrices varies
between 20\% and 100\% and the discount factor is set to 0.9, 0.98, and 0.995. This
results in 27 problem instances in each example, which are then solved using the
standard value iteration algorithm and its three variants with and without our
acceleration operators.

\begin{exl} \label{ex:dense}
Consider MDPs with 500 states and up to 100 actions per state. The transition
probability matrices were generated using a uniform random number generator and
non-zero elements were uniformly placed in the matrix. The density of non-zero
elements of the matrices varies from 20\% to 100\%. The discount factor is set to
0.9, 0.98, and 0.995 and the tolerance used in the stopping condition is
$\eps=10^{-3}$.
\end{exl}

The computational results of Example 1 are presented in Table~\ref{tab:1}, where the
first and the second numbers in each cell are the number of iterations and the CPU
time to converge to fixed points, respectively. PAVI improves the computational
efficiency measured in the CPU time up to 605 times (PAVI combined with Jacobi value
iteration when the density is 50\% and the discount factor is 0.995). In general,
PAVI has a larger improvement than LAVI, although LAVI is better in terms of the worst
case comparison. Both operators improve the computation as the discount factor
approaches 1, which is the case where the value iteration algorithm and its variants
suffer from slow convergence.

\vspace{1in}
\begin{center}
Table~\ref{tab:1} goes here.
\end{center}
\vspace{1in}

\begin{exl} \label{ex:diag}
Consider MDPs with 500 states and up to 100 actions per state.The transition
probability matrices are band matrix - matrix with non-zero elements only in a band
centering around the diagonal of the matrix - with a bandwidth varying from 100 to
450 (which gives us the density of the matrix varying from 20\% to 90\%). When the
density is set to 100\%, Example 2 becomes identical to Example 1.  The transition
probabilities inside the band were generated using a uniform random number generator.
The discount factor is set to 0.9, 0.98, and 0.995, and the tolerance used in
the stopping condition is $\eps=10^{-3}$.
\end{exl}

The computational results are presented in Table~\ref{tab:2}. The general trend in
Example 2 is similar to that in the earlier example except for the acceleration being
less dramatic. The maximum acceleration in terms of the CPU time is observed when the
discount factor is 0.995 and the density is 80\%, and PAVI accelerates the convergence
217 times over Jacobi value iteration. Similar to Example 1, the PAVI outperforms
LAVI, both operators perform best when combined with the standard value iteration,
and both accelerators perform better as the discount factor approaches 1 and the
density approaches 100\%.

Let us now present the brief analysis of the above numerical results.
Based on Example 1 and Example 2, we can conclude that PAVI generally performs better
than LAVI, and that the acceleration from PAVI is especially beneficial when it is
combined with the standard value iteration and Jacobi methods. LAVI exhibits a
relatively minor acceleration with these two methods and shows stronger performance
when it is combined with Gauss-Seidel and Gauss-Seidel-Jacobi algorithm. We propose the the following explanation to this behavior. As one can notice, the impact of the projection operator is more significant when the Markovian operator moves the points deeper inside the set V rather then closer to fixed point. Because of this, combining the projection operator with the standard VI and Jacobi VI is more advantageous than combining it with Gauss-Seidel and Gauss-Seidel-Jacobi VI.
Moreover, it is assumed that than the higher the density of the transition probability matrix the deeper inside of the set $V$ the point will be moved by VI in any of the above forms. This assumption is related to the result of Lemma \ref{lem:fully dense}, where fully dense matrices are the limiting case of the dense matrices. One can see, especially from results in Table \ref{tab:2}, the denser the transition probability matrix, the better the performance of PAVI combined with standard and Jacobi VI. One more remark should be made.

\begin{rem}
  \label{r:PAGSJ=GSJ}
  Notice that in some cases it turns out that PAVI+GSJ requires significantly more iterations than PAVI+GS. The reason is that the coordinate
  $\abs{S}$ of $u=\TGSJ{}v$ is expressed in terms of $u_j$,
  $j<\abs{S}$. Therefore for a certain $a\in{}A(\abs{S})$ we have the
  equality
  $$
  u_{\abs{S}}= \Bigl[ r(i,a)+\lm\sum_{j<\abs{S}} p_{\abs{S},j}(a) u_j \Bigr ]/ {\Bigl [ 1-\lm p_{|S||S|}(a)\Bigr ]},
  $$
  which easily implies that $u\in\partial{}V$. Therefore, though it is not
  guaranteed, it is possible that
  $P\TGSJ{}v=\TGSJ{}v$ for all $v\in{}V$.
\end{rem}

On the other hand, one can notice that Gauss-Seidel and Gauss-Seidel-Jacobi VI usually converge faster than the standard VI and Jacobi VI, but they iterate points closer to facets of set V, which reduces the impact of using projection operators. As suggested by geometry of set V and supported by the numerical examples, combining Gauss-Seidel and Gauss-Seidel-Jacobi VI with a linear extension operator is more advantageous than combining standard VI and Jacobi VI with it.

For almost all of the cases in Table 1, the combinations of PAVI + S and PAVI + J are the best combinations; however, for Table 2 one can see that for cases with sparse transition probability matrices and  $\lm$ close to 1, the combinations LAVI + GS and LAVI + GSJ are more advantageous. The numerical studies show that combinations of PAVI + S and PAVI + J perform better when the transition probability matrices are relatively dense ( $10\%$ or more), while for sparse matrices the combinations LAVI + GS and LAVI + GSJ may perform better. For very sparse matrices ($0.1\%$ or less) the performance of all variants of GAVI is relatively poor.

In general, the acceleration tends to be more dramatic as the density increases and as the discount
factor approaches 1. This is encouraging, since the standard value iteration
algorithm and its variants have been known to suffer from slow convergence when the
discount factor is close to 1 and requires more computation as the transition
probability matrix becomes denser. Therefore, our accelerators will perform
excellently when they are needed the most.

\subsection{Computational Complexity and Savings}
We now present the computational cost of the proposed operators. For the standard VI,
when the transition probability matrices are fully dense, each iteration will take
$C|S|^2$ (where $C$ is the average number of actions per state) multiplications and
divisions. With sparse transition probability matrices, this number can be estimated
as $NC|S|$ (where $N$ is the average number of non-zero entries per row of the
transition probability matrices).

The additional effort required in GAVI is due to the acceleration operator used in
{\bf Step 1} of GAVI. Both accelerators require the computation of
$\sum_{j}p_{ij}(a)Tw^{n-1}_j $ for all $i\in{}I$ for all $a\in{}A(i)$ in $n$-th
iteration of GAVI. Notice that we will have to know $\sum_{j}p_{ij}(a)w^{n}_j$ in the
next iteration in order to compute $Tw^{n}$. For example, for
PAVI $w^n=\alpha^* Tw^{n-1}$ and $ \sum_{j}p_{ij}(a)w^{n}_j =\alpha^*
\sum_{j}p_{ij}(a)Tw^{n-1}_j$. Therefore, we can store $\sum_{j}p_{ij}(a)Tw^{n-1}_j$
between iterations to avoid redundant computation. A similar approach can be applied
to the three variants of acceleration operators studied in this paper. Therefore,
each iteration of GAVI requires $C|S||S+1|$ multiplication and division, which is
just slightly more than that of standard value iteration.

From the above we can draw the following conclusion: the additional computation due to the accelerating operators is
marginal. Though it varies among different problem instances, the average additional
computation per iteration observed in 54 instances from Example 1 and Example 2, each
of which is solved using 12 different approaches, is 12\% with PAVI and 15\% with
LAVI.

\section{The Expected Total Reward and Other MDPs} \label{sec:total-reward}
In this section we will show that GAVI discussed in Section \ref{sec:accel-value} can
also be used with other types of MDP: MDPs with the expected total reward criterion,
semi-Markov decision processes, and continuous-time MDP.

Consider an infinite horizon MDP with expected total reward criterion. One can notice
that both the coordinate translation and the rewards adjustment described in
Remark~\ref{r:non-negativity} cannot be applied to this model. This is because of
$\lm=1$ for the coordination translation and the non-zero reward at absorbing states for
the reward adjustment. Therefore, the proposed techniques are applicable only to
positive models (and also to negative models with trivial changes in the definition
of the acceleration operators).

Recall that the value iteration algorithm for MDP with expected total reward
criterion is based on successive application of an operator $L$ to a vector
$v\in\R^{|S|}$~\cite{Puterman:MDP94},
\begin{equation} \label{eq:3aTR}
(Lv)_i=\max_{a\in{A(i)}}\left\{r(i,a)+\sum_j p_{ij}(a)v_j\right\}.
\end{equation}
Since the operator $L$ in Equation (\ref{eq:3aTR}) is not a contraction mapping, the
convergence of successive approximation algorithms using the operator $L$ is not
guaranteed. However, with $v^0$ satisfying the conditions $0 \le v^0 \le v^*$ or $v*
\le v^0$, the sequence $\{v^n\}_{n=0}^\infty$ converges monotonically to $v^*$
\cite{Puterman:MDP94}.

The LP formulation equivalent to an MDP with expected total reward criterion is given
as follows:
\begin{equation} \label{eq:LP4Total}
\min \left\{\sum_i v_i \mid{} v_i-\sum_jp_{ij}(a)v_j\ge{}r(i,a), \forall i\in{S},
\forall a\in A(i), v \in \R^{|S|} \right\},
\end{equation}
where $v^*$, the optimal solution of the LP, yields the fixed point of the MDP.
Analogous with the discounted MDP considered in Section \ref{sec:intro}, the feasible
set $W$ of the LP in Equation (\ref{eq:LP4Total}) can be described as follows:
\begin{equation}
  W=\{v\in\R^{|S|} \mid v\ge{}Lv\}. \label{eq:VJTR}
\end{equation}
The following Lemma states that the same approaches as introduced in Section
\ref{sec:accel-value} are applicable to MDPs with the total expected reward.
\begin{lem} \label{lem:invariance4total}\
$W$ is invariant under $L$. That is, $LW \subset W$.
\end{lem}
With the invariance of $W$ under $L$, not only GAVI(General Accelerated Value
Iteration) algorithm but also all of its variants can be applied directly to MDPs
with the expected total reward criterion.

\begin{rem} \label{r:Variants of VI}
Jacobi and Gauss-Seidel-Jacobi variants of the value iteration algorithm are no
longer valid for MDPs with the expected total reward criterion.
\end{rem}

The value iteration algorithm for MDPs with the expected total reward has been known
for its very slow convergence. Our acceleration operators, which improve the
efficiency by an increasingly margin as the discount factor approaches 1, are an
attractive option for MDPs with the expected total reward, which can be seen as MDPs
with the discounted reward with a discount factor 1. Example \ref{ex:AV} below
confirms this speculation.

\begin{exl} \label{ex:AV} Consider an MDP with 5 states and up to 5 actions per state.
Let one of the states be an absorbing state with zero reward. The standard value
iteration takes 293 iterations, while our acceleration value iteration with
Projective Operator takes only 15 iterations before meeting the stopping condition
($\eps=10^{-3}$).
\end{exl}

Now we will make a short remark about semi-Markov and continuous-time decision
process. It is widely known that discounted infinite-horizon semi-Markov decision
processes and continuous-time MDPs can be transformed into models similar to a discrete
time discounted model with rewards $r(i, a)$ depending only on state at decision
epoch and action~\cite{Lippman:OR75p687,Puterman:MDP94}. Therefore, our acceleration
operations can be used for these types of MDPs after they are converted into discrete
time MDPs with the total expected discounted reward.

\section{Conclusions} \label{sec:conclusions}
Inspired by the monotone behavior of the contraction mapping operator used in the
value iteration algorithm within the feasible set of the linear programming problem
equivalent to the given MDP, we propose a class of operators that can be used in
combination with the standard contraction mapping and its variants, such as Jacobi,
Gauss-Seidel, and Gauss-Seidel-Jacobi methods, to improve the computational
efficiency. Two acceleration operators, Projective Operator and Linear Extension
Operator, are particularly proposed and combined into the standard value iteration
algorithm and its variants. The savings due to the acceleration have been dramatic and
the maximum savings is up to 607 time faster than the case without our accelerating
operator. The numerical studies show that combinations of PAVI + S and PAVI + J perform better when the transition probability matrices are relatively dense ( $10\%$ or more), while for sparse matrices the combinations LAVI + GS and LAVI + GSJ may perform better.
It is especially interesting to mention that the savings become more
dramatic when encountering problems for which the standard value iteration algorithm suffers
from slow convergence (when the discount factor approaches 1 and when the transition
probability matrix becomes dense). The strong performance with the discount factor
close to 1 can be of special interest to those who are investigating Blackwell
Optimality, which assures that a Blackwell optimal policy is optimal for all discount
factors $\lambda_b < 1$ close enough to
1~\cite{Blackwell:AMS62p719,Veinott:AMS69p1635}.

\appendix

\section{Proofs} \label{sec:proofs}

\begin{lem}  \label{l:monotonicity}
\noindent \begin{enumerate}
  \renewcommand{\labelenumi}{(\roman{enumi})}
  \item  If $u\ge{}v$, then
        $\;\; Tu\ge{}Tv$,  $T_J u\ge{}T_J v$,  $\TGS u\ge{}\TGS v$,  and $\TGSJ u\ge{}\TGSJ v$.
  \item  If $u>{}v$, then
        $\;\; Tu>{}Tv$, $T_J u>{}T_J v$, $\TGS u>{}\TGS v$, and $\TGSJ u>{}\TGSJ v$.
\end{enumerate}
\end{lem}
\begin{proof}[Proof of Lemma~\ref{l:monotonicity}]
Let us first prove that $u\ge{}v$ implies $Tu\ge{}Tv$.
\[
Tu_i=\max_{a\in A(i)}\left\{r(i,a)+\lm\sum_{j}p_{ij}(a)u_j\right\} \ge \max_{a\in
A(i)}\left\{r(i,a)+\lm\sum_{j}p_{ij}(a)v_j\right\} = Tv_i \text{ for all } i\in S.
\]
For part (ii) the inequality can be simply replaced with a strict inequality. The
same proof applies for $T_Ju\ge T_Jv$. Let $w=\TGS u$ and $\xi=\TGS v$. Then $ w_1=
\TGS u_1 = Tu_1 \ge Tv_1 = \xi_1 \text{ for all } i\in S. $ By induction, assuming
that $w_k \ge \xi_k$ for all $k<i$, we get for
  $j=i$
  \begin{align*}
    w_i&=\max_{a\in A(i)}\biggl(r(i,a)+\lm\sum_{j<i}p_{ij}(a)w_j +
    \lm\sum_{j \ge i}p_{ij}(a)u_j \biggr)\\
         &\ge\max_{a\in A(i)}\biggl(r(i,a)+\lm\sum_{j<i}p_{ij}(a)\xi_j
         + \lm\sum_{j \ge i}p_{ij}(a)v_j \biggr)=(\TGS v)_i=\xi_i.
  \end{align*}
This proof can be trivially modified for $\TGSJ u\ge{}\TGSJ v, \TGS u>{}\TGS v$, and
$\TGSJ u>{}\TGSJ v$.
\end{proof}

\begin{rem}
Notice that proof of Lemma~\ref{l:monotonicity} is good even when $\lm=1$.
\end{rem}

\begin{proof}[Proof of Lemma~\ref{lem:invariance}]
Let $v \in V$ and $u=Tv$. By definition of set $V$, $u = Tv \le v$. By monotonicity
shown in Lemma~\ref{l:monotonicity}, $Tu \le Tv=u$. Thus, $u \in V$.
\end{proof}

\begin{proof}[Proof of Theorem~\ref{l:Pu_satisfy12}]
Condition (A) is satisfied trivially since $\al v \in V$ for any $v\in V$ by the
definition of $Z$ given in (\ref{eq:1DLP}). Now we have to show that $Z$ satisfies
condition (B). We know $\al=1$ is feasible to the linear program (\ref{eq:1DLP})
since $v\in V$ (or $Tv\le v$). Since $\sum_i v_i\ge 0$ due to $r(i,a)\ge 0,\; \forall
i, a$, we have $\al^* \le 1$. Therefore, $Zv = \al^* v \le v$.
\end{proof}

\begin{proof} [Proof of Lemma~\ref{lem:fully dense}]
Suppose that the transition probability matrix is fully dense (i.e., $p_{ij}(a)>0, \;
\forall i,j,a$). Then if $v\ge u=Tv$ and $v \neq v^*$, where $v^*$ is the fixed point
of $T$, there exists $k$ such that $v_k
> u_k$. As a result,
  \begin{align*}
    Tu_i&=\max_{a\in A(i)}\left\{r(i,a)+\lm\sum_{j\neq k}p_{ij}(a)u_j
    +\lm p_{ik}(a)u_{k}\right\}\\
    &<\max_{a\in A(i)}\left\{r(i,a)+
    \lm\sum_{j\neq k}p_{ij}(a)v_j +\lm p_{ik}(a)v_{k}\right\} = v_i
    \text{ for all } i\in S.
 \end{align*}
\end{proof}

\begin{proof}[Proof of Lemma~\ref{lem:interior}]
Let $v \in int(V)$ and $u=Tv$. By definition of $int(V)$,  $u = Tv < v$.  By
monotonicity shown in Lemma~\ref{l:monotonicity}, $Tu < Tv=u$. Thus, $u \in int(V)$.
\end{proof}

\begin{proof}[Proof of Theorem~\ref{thm:L_w_satisfy12}]
For $v\in V$, $Z(v)=v+\al^*(Tv-v)$, where $\al^*$ is an optimal solution to the
linear program in \ref{eq:1DLP2}. Since $v+\al^*(Tv-v)$ is feasible to the linear
program, we have $T(v+\al^*(Tv-v))\le v+\al^*(Tv-v)$. This suffices Condition (A). By
$Tv=u\in V$, $\al=1$ is feasible. Since $Tv\le v$ and $\al=1$ is feasible, $\al^*>0$.
Hence, $Zv \le v$ and Condition (B) is satisfied.
\end{proof}

\begin{proof}[Proof of Lemma~\ref{lem:invariance of splittings}]
This proof is similar to the proofs of Lemma~\ref{lem:invariance} and ~\ref{lem:interior}.
\end{proof}

\begin{proof} [Proof of Lemma~\ref{l:relationships_between_V's}]
First let's show that $V_J \in V$. From (\ref{eq:TJ}), $v^{n+1}=T_J v^n$, where
\begin{align*}
v^{n+1}_i &=\max_{a \in A(i)} \left\{\left[r(i,a)+\lambda \sum_{j\neq
            i}p_{ij}(a)v^n_j\right] /\left[1-\lambda p_{ii}(a)\right]\right\} \\
       &\ge \left[r(i,a)+\lambda \sum_{j\neq i}p_{ij}(a)v^n_j\right]
               /\left[1-\lambda p_{ii}(a)\right], \forall a \in A(i),
\end{align*}
which then implies
\[
v^{n+1}_i-v^{n+1}_i\lambda p_{ii}(a)\ge r(i,a)+\lambda \sum_{j\neq i}p_{ij}(a)v^n_j
\]
or
\[
v^{n+1}_i \ge r(i,a)+\lambda \sum_{j}p_{ij}(a)v^n_j.
\]
Therefore, $V_J \subset V$. Repeating the same steps and moving backward it is easy
to show that $V \subset V_J$. Similarly we can show $V_{GS}=\VGSJ$. In order to prove
inclusion $V \subset \VGS$, it is sufficient to show that if $v\ge Tv$, then $v\ge
T_{GS}v$. For $v\in V$, let $u=\TGS{v}$ and $w=Tv$. Then, $w_j\le{}v_j$ for all $j$
and $u_1=w_1$. Assume that $u_k\le{}w_k$ for all $k<i$, then
\begin{align*}
 u_i&=\max_{a\in A(i)}\biggl(r(i,a)+\lm\sum_{j<i}p_{ij}(a)u_j +
      \lm\sum_{j \ge i}p_{ij}(a)v_j \biggr)\\
    &\le\max_{a\in{}A(i)}\biggl(r(i,a)+\lm\sum_{j<i}p_{ij}(a)w_j +
      \lm\sum_{j \ge i}p_{ij}(a)v_j \biggr)\\
    &\le\max_{a\in{}A(i)}\biggl(r(i,a)+\lm\sum_{j<i}p_{ij}(a)v_j +
      \lm\sum_{j \ge i}p_{ij}(a)v_j \biggr)=(Tv)_i=w_i.
\end{align*}
By induction, $w\ge u$, implying $v\ge Tv=w\ge u = T_{GS}v$.
\end{proof}

\begin{proof}[Proof of Lemma~\ref{l:TX(VX)subsetV}]
For $v\in\VGS$, let $u=\TGS{v}$.  By Lemma~\ref{l:monotonicity}, $v\ge\TGS{}v=u$. By
replacing ``$\max_{a\in A(i)}$'' by inequalities, similar to the argument used in the
proof of Lemma~\ref{l:relationships_between_V's}, and by $v\ge u$, we have
\begin{align*}
 u_i &\ge{}r(i,a)+\lm\sum_{j<i}p_{ij}(a)u_j+\lm\sum_{j \geq i}p_{ij}(a)v_j
   \text{ for all $i\in{}I$ for all $a\in{}A(i)$} \\
     &\ge r(i,a)+\lm\sum_{j<i}p_{ij}(a)u_j+
  \lm\sum_{j\ge{}i}p_{ij}(a)u_j \text{ for all $i\in{}I$ for all $a\in{}A(i)$,}
\end{align*}
which is equivalent to $u\ge{}Tu$, or $u\in{V}$. If $v\in\VGSJ$, then letting
$u=\TGSJ$ by similar argument we get $u\ge{}T_Ju$. Since $V_J=V$, $u\in{V}$.
\end{proof}

\begin{proof}[Proof of Theorem~\ref{thm:V_invariance_under_TJ,TGS,TGSJ}]
By monotonicity shown in Lemma~\ref{l:monotonicity}, $T_JV_J \subset V_J$. Therefore,
by $V=V_J$ in Lemma~\ref{l:relationships_between_V's}, $T_JV \subset V$. By
Lemma~\ref{l:relationships_between_V's} and Lemma~\ref{l:TX(VX)subsetV},
\[
T_{GS}V \subset T_{GS}V_{GS} \subset V \;\;\text{and}\;\; T_{GSJ}V\subset
T_{GSJ}V_{GSJ}\subset V.
\]
\end{proof}

\begin{proof}[Proof of Lemma~\ref{lem:invariance4total}]
The proof is similar to proofs of Lemma~\ref{lem:invariance} and ~\ref{lem:interior}.
\end{proof}

\section*{Acknowledgement}
We express our deepest condolences to the late Dr. Dmitry Khmelev, who was a
postdoctoral fellow in the Department of Mathematics at University of Toronto before
he passed away in 2004.

\bibliographystyle{abbrv}
%\bibliographystyle{named}
%\bibliography{ref}

\begin{center}
\begin{table} [ht]
\fontsize{9}{7} \selectfont
\begin{tabular}{|c|c|c|c|c|c|c|c|c|c|c|c|} \hline
 \multicolumn{2}{|c|}{$\lm$}&\multicolumn{3}{|c|}{0.9}&\multicolumn{3}{|c|}{0.98}&\multicolumn{3}{|c|}{0.995}\\\hline
 Dens.$^\dag$&{VI}&None$^\ddag$&PAVI$^\ddag$&LAVI$^\ddag$&None$^\ddag$&PAVI$^\ddag$&LAVI$^\ddag$&None$^\ddag$&PAVI$^\ddag$&LAVI$^\ddag$\\\hline
 & S & 201 / 49& 7 / 2& 57 / 15& 1164 / 287& 8 / 2& 339 / 91& 5104 / 1257& 8 / 2& 1504 / 397\\\cline{2-11}
& J & 201 / 56  & 7 / 2& 56 / 17& 1163 / 324& 8 / 3& 339 / 102& 5097 / 1427& 8 / 2& 1508 / 451 \\\cline{2-11}
 \raisebox{1.9ex}[0cm][0cm]{100} & GS &  75 / 18& 58 / 26& 33 / 15& 419 / 103& 317 / 139& 180 / 80& 1818 / 449& 1188 / 522& 774 / 343\\\cline{2-11}
&GSJ & 74 / 20& 59 / 26& 28 / 13& 417 / 117& 336 / 148& 110 / 50& 1811 / 503& 1485 / 653& 379 / 169\\\hline
 & S & 201 / 47& 7 / 1& 58 / 15& 1164 / 273& 7 / 2& 335 / 84& 5104 / 1196& 8 / 2& 1498 / 377\\\cline{2-11}
& J & 201 / 53  & 6 / 2& 57 / 17& 1163 / 308& 7 / 2& 336 / 95& 5098 / 1355& 8 / 2& 1493 / 418\\\cline{2-11}
 \raisebox{1.9ex}[0cm][0cm]{90} & GS &  75 / 17& 57 / 24& 26 / 11& 418 / 98& 306 / 127& 163 / 68& 1817 / 427& 1081 / 447& 743 / 310\\\cline{2-11}
& GSJ & 74 / 20& 58 / 24& 20 / 9& 417 / 111& 333 / 139& 97 / 42& 1810 / 478& 1474 / 614& 310 / 131\\\hline
 & S & 201 / 43& 7 / 2& 58 / 13& 1164 / 248& 8 / 1& 342 / 78& 5102 / 1081& 9 / 2& 1521 / 345\\\cline{2-11}
& J & 201 / 48  & 7 / 2& 57 / 15& 1162 / 277& 8 / 2& 340 / 89& 5095 / 1216& 9 / 2& 1504 / 393\\\cline{2-11}
 \raisebox{1.9ex}[0cm][0cm]{80} & GS&  75 / 16& 56 / 21& 32 / 13& 418 / 89& 275 / 101& 180 / 68& 1815 / 385& 788 / 290& 768 / 285\\\cline{2-11}
& GSJ & 74 / 17& 58 / 22& 26 / 10& 416 / 98& 334 / 124& 121 / 46& 1808 / 429& 1477 / 546& 399 / 149\\\hline
 & S & 201 / 37& 7 / 1& 59 / 12& 1163 / 223& 8 / 1& 348 / 70& 5100 / 973& 9 / 2& 1530 / 310\\\cline{2-11}
& J & 201 / 42  & 7 / 1& 58 / 14& 1162 / 246& 8 / 1& 345 / 77& 5093 / 1080& 9 / 2& 1538 / 354\\\cline{2-11}
 \raisebox{1.9ex}[0cm][0cm]{70} & GS&  75 / 14& 59 / 19& 31 / 11& 418 / 81& 330 / 108& 174 / 57& 1815 / 342& 1365 / 447& 769 / 253\\\cline{2-11}
& GSJ & 74 / 16& 59 / 19& 31 / 11& 416 / 88& 337 / 111& 142 / 47& 1807 / 382& 1488 / 489& 443 / 147\\\hline
 & S & 201 / 33& 8 / 1& 62 / 11& 1163 / 190& 9 / 1& 356 / 62& 5098 / 835& 10 / 1& 1558 / 272\\\cline{2-11}
& J & 201 / 37  & 8 / 1& 61 / 12& 1161 / 211& 9 / 2& 353 / 70& 5092 / 930& 10 / 2& 1549 / 304\\\cline{2-11}
 \raisebox{1.9ex}[0cm][0cm]{60} & GS&  74 / 12& 60 / 17& 33 / 9& 417 / 71& 342 / 95& 108 / 31& 1813 / 297& 1510 / 419& 418 / 117\\\cline{2-11}
& GSJ & 74 / 13& 60 / 17& 33 / 9& 416 / 75& 341 / 96& 106 / 30& 1806 / 327& 1505 / 421& 441 / 125\\\hline
 & S & 201 / 31& 8 / 1& 62 / 10& 1165 / 180& 9 / 2& 353 / 58& 5106 / 789& 9 / 2& 1546 / 257\\\cline{2-11}
& J & 201 / 34  & 8 / 1& 61 / 11& 1163 / 199& 9 / 2& 353 / 66& 5098 / 873& 9 / 1& 1543 / 292\\\cline{2-11}
 \raisebox{1.9ex}[0cm][0cm]{50} & GS&  75 / 12& 57 / 15& 32 / 9& 419 / 65& 285 / 74& 180 / 48& 1820 / 281& 855 / 222& 788 / 208\\\cline{2-11}
& GSJ & 74 / 13& 59 / 15& 32 / 9& 417 / 72& 336 / 87& 108 / 29& 1812 / 312& 1486 / 387& 322 / 86\\\hline
 & S & 201 / 28& 8 / 1& 62 / 9& 1164 / 163& 9 / 1& 355 / 52& 5105 / 715& 10 / 1& 1571 / 235\\\cline{2-11}
& J & 201 / 31  & 8 / 1& 61 / 9& 1163 / 178& 9 / 2& 353 / 58& 5099 / 784& 10 / 2& 1565 / 262\\\cline{2-11}
 \raisebox{1.9ex}[0cm][0cm]{40} & GS&  75 / 10& 59 / 13& 32 / 7& 418 / 61& 336 / 77& 101 / 24& 1816 / 255& 1486 / 340& 111 / 26\\\cline{2-11}
& GSJ & 74 / 12& 58 / 13& 30 / 7& 417 / 64& 335 / 78& 105 / 24& 1810 / 277& 1481 / 342& 139 / 33\\\hline
 & S & 201 / 25& 8 / 1& 64 / 8& 1164 / 145& 9 / 2& 358 / 48& 5102 / 635& 10 / 1& 1584 / 210\\\cline{2-11}
& J & 201 / 27  & 8 / 1& 63 / 8& 1162 / 157& 9 / 1& 356 / 55& 5096 / 689& 10 / 1& 1580 / 221\\\cline{2-11}
 \raisebox{1.9ex}[0cm][0cm]{30} & GS&  74 / 9& 58 / 12& 32 / 7& 418 / 52& 334 / 66& 82 / 16& 1813 / 232& 1479 / 291& 113 / 23\\\cline{2-11}
& GSJ & 74 / 10& 58 / 12& 32 / 7& 416 / 56& 333 / 66& 71 / 15& 1807 / 245& 1473 / 292& 134 / 27\\\hline
 & S & 201 / 21& 8 / 1& 64 / 7& 1163 / 122& 9 / 1& 361 / 40& 5101 / 534& 10 / 1& 1605 / 178\\\cline{2-11}
& J & 201 / 22  & 8 / 1& 63 / 8& 1162 / 130& 9 / 1& 360 / 47& 5093 / 574& 10 / 1& 1600 / 199\\\cline{2-11}
 \raisebox{1.9ex}[0cm][0cm]{20} & GS&  74 / 8& 58 / 9& 31 / 5& 417 / 44& 332 / 53& 86 / 14& 1813 / 190& 1471 / 233& 114 / 19\\\cline{2-11}
& GSJ & 74 / 9& 58 / 10& 23 / 4& 415 / 46& 331 / 54& 70 / 11& 1804 / 201& 1463 / 236& 190 / 31\\\hline
\end{tabular}

\vspace{0.1in} {\footnotesize $\dag$: The density of the transition probability matrix (\%).\\
$\ddag$: Two values in these columns correspond to the
number of iterations and CPU time of the corresponding algorithm.} \vspace{-0.1in}
\caption{The number of iterations and CPU time of the standard, Jacobi, Gauss-Seidel,
and Gauss-Seidel-Jacobi value iteration algorithms with and without an accelerating
operator applied to a family of MDPs from Example~\ref{ex:dense}} \label{tab:1}
\end{table}
\end{center}

\begin{center}
\begin{table} [ht]
\fontsize{9}{7} \selectfont
\begin{tabular}{|c|c|c|c|c|c|c|c|c|c|c|c|} \hline
 \multicolumn{2}{|c|}{$\lm$}&\multicolumn{3}{|c|}{0.9}&\multicolumn{3}{|c|}{0.98}&\multicolumn{3}{|c|}{0.995}\\\hline
Density$^\dag$&{VI}&None$^\ddag$&PAVI$^\ddag$&LAVI$^\ddag$&None$^\ddag$&PAVI$^\ddag$&LAVI$^\ddag$&None$^\ddag$&PAVI$^\ddag$&LAVI$^\ddag$\\\hline
 & S & 201 / 36& 15 / 3& 65 / 12& 1164 / 209& 20 / 4& 385 / 75& 5102 / 919& 23 / 4& 1711 / 336\\\cline{2-11}
& J & 201 / 42  & 15 / 3& 64 / 13& 1162 / 244& 20 / 4& 384 / 79& 5092 / 1081& 23 / 5& 1710 / 343\\\cline{2-11}
 \raisebox{1.9ex}[0cm][0cm]{90} & GS &  75 / 14& 59 / 16& 29 / 8& 419 / 76& 334 / 89& 171 / 46& 1818 / 327& 1465 / 387& 752 / 201\\\cline{2-11}
& GSJ & 75 / 16& 58 / 16& 21 / 6& 417 / 87& 333 / 89& 123 / 34& 1815 / 382& 1473 / 391& 321 / 87\\\hline
 & S & 201 / 33& 15 / 2& 64 / 10& 1164 / 179& 22 / 3& 378 / 63& 5103 / 785& 26 / 4& 1678 / 282\\\cline{2-11}
& J & 201 / 36  & 15 / 2& 63 / 11& 1162 / 209& 22 / 4& 377 / 69& 5093 / 912& 26 / 4& 1679 / 298\\\cline{2-11}
 \raisebox{1.9ex}[0cm][0cm]{80} & GS&  75 / 12& 56 / 13& 33 / 8& 419 / 65& 252 / 58& 173 / 40& 1820 / 280& 621 / 143& 777 / 181\\\cline{2-11}
& GSJ & 75 / 13& 59 / 14& 29 / 7& 417 / 74& 335 / 78& 108 / 27& 1816 / 328& 1478 / 343& 297 / 70\\\hline
 & S & 202 / 26& 17 / 3& 66 / 9& 1164 / 152& 26 / 4& 392 / 56& 5101 / 704& 32 / 5& 1746 / 247\\\cline{2-11}
& J & 201 / 30  & 17 / 3& 65 / 9& 1161 / 175& 26 / 4& 391 / 58& 5090 / 766& 32 / 4& 1745 / 251\\\cline{2-11}
 \raisebox{1.9ex}[0cm][0cm]{70} & GS&  75 / 9& 58 / 11& 32 / 7& 420 / 55& 278 / 56& 175 / 36& 1820 / 235& 786 / 156& 766 / 154\\\cline{2-11}
& GSJ & 75 / 11& 59 / 12& 32 / 7& 417 / 63& 337 / 68& 145 / 30& 1812 / 274& 1487 / 300& 367 / 78\\\hline
 & S & 201 / 23& 29 / 3& 68 / 8& 1163 / 124& 47 / 6& 398 / 47& 5096 / 544& 57 / 6& 1766 / 207\\\cline{2-11}
& J & 201 / 25  & 29 / 3& 67 / 8& 1160 / 144& 47 / 5& 399 / 48& 5083 / 632& 57 / 7& 1767 / 209\\\cline{2-11}
 \raisebox{1.9ex}[0cm][0cm]{60} & GS&  75 / 8& 56 / 9& 32 / 6& 419 / 44& 226 / 37& 182 / 30& 1817 / 193& 477 / 78& 785 / 131\\\cline{2-11}
& GSJ & 75 / 10& 60 / 10& 34 / 5& 416 / 52& 341 / 57& 127 / 21& 1809 / 225& 1503 / 248& 361 / 61\\\hline
 & S & 202 / 19& 45 / 4& 68 / 7& 1163 / 110& 80 / 8& 396 / 42& 5096 / 481& 99 / 10& 1756 / 185\\\cline{2-11}
& J & 201 / 23  & 45 / 5& 67 / 9& 1160 / 127& 80 / 9& 396 / 49& 5082 / 556& 99 / 10& 1759 / 210\\\cline{2-11}
 \raisebox{1.9ex}[0cm][0cm]{50} & GS&  76 / 7& 57 / 9& 31 / 5& 419 / 40& 251 / 39& 162 / 25& 1818 / 171& 605 / 90& 725 / 111\\\cline{2-11}
& GSJ & 75 / 8& 59 / 9& 32 / 6& 416 / 46& 335 / 50& 124 / 20& 1806 / 207& 1479 / 220& 239 / 37\\\hline
 & S & 202 / 16& 51 / 4& 69 / 7& 1163 / 93& 108 / 9& 400 / 35& 5095 / 429& 143 / 12& 1765 / 155\\\cline{2-11}
& J & 201 / 19  & 50 / 5& 68 / 9& 1159 / 106& 108 / 10& 399 / 40& 5077 / 483& 142 / 12& 1767 / 173\\\cline{2-11}
 \raisebox{1.9ex}[0cm][0cm]{40} & GS&  76 / 6& 56 / 7& 32 / 4& 420 / 33& 236 / 29& 179 / 23& 1819 / 143& 526 / 65& 742 / 94\\\cline{2-11}
& GSJ & 75 / 7& 60 / 8& 32 / 4& 416 / 40& 334 / 41& 141 / 19& 1806 / 168& 1472 / 184& 235 / 30\\\hline
 & S & 203 / 13& 62 / 4& 69 / 5& 1163 / 73& 170 / 12& 399 / 28& 5094 / 318& 252 / 17& 1763 / 123\\\cline{2-11}
& J & 202 / 14  & 62 / 4& 68 / 6& 1158 / 83& 169 / 13& 398 / 32& 5071 / 365& 251 / 18& 1765 / 136\\\cline{2-11}
 \raisebox{1.9ex}[0cm][0cm]{30} & GS&  77 / 5& 57 / 6& 35 / 3& 421 / 26& 235 / 24& 174 / 18& 1823 / 114& 514 / 51& 769 / 79\\\cline{2-11}
& GSJ & 76 / 5& 60 / 6& 35 / 4& 416 / 31& 333 / 34& 159 / 16& 1801 / 129& 1466 / 148&  408 / 21 \\\hline
 & S & 204 / 9& 76 / 3& 68 / 3& 1163 / 52& 283 / 13& 400 / 21& 5085 / 226& 525 / 25& 1772 / 89\\\cline{2-11}
& J & 202 / 10  & 75 / 4& 68 / 3& 1155 / 59& 280 / 15& 399 / 23& 5052 / 258& 520 / 27& 1767 / 101\\\cline{2-11}
 \raisebox{1.9ex}[0cm][0cm]{20} & GS&  78 / 4& 57 / 4& 36 / 3& 422 / 18& 216 / 16& 182 / 13& 1823 / 81& 420 / 31& 768 / 59\\\cline{2-11}
& GSJ & 77 / 4& 60 / 4& 34 / 2& 414 / 22& 330 / 25& 149 / 12& 1789 / 93& 1448 / 107& 298 / 23\\\hline
 & S & 204 / 5& 82 / 3& 69 / 2& 1162 / 32& 361 / 11& 400 / 13& 5065 / 138& 1091 / 31& 1772 / 56\\\cline{2-11}
& J & 201 / 6  & 80 / 3& 68 / 2& 1147 / 34& 355 / 13& 398 / 15& 4997 / 151& 1071 / 37& 1767 / 63\\\cline{2-11}
 \raisebox{1.9ex}[0cm][0cm]{10} & GS&  80 / 3& 57 / 2& 37 / 2& 430 / 11& 248 / 12& 189 / 10& 1838 / 49& 674 / 32& 796 / 40\\\cline{2-11}
& GSJ & 78 / 3& 61 / 3& 37 / 2& 414 / 13& 332 / 16& 164 / 8& 1770 / 54& 1438 / 69& 429 / 22\\\hline
 \end{tabular}

\vspace{0.1in} {\footnotesize $\dag$: The density of the transition probability matrix (\%).\\
$\ddag$: Two values in these columns correspond to the
number of iterations and CPU time of the corresponding algorithm.} \vspace{-0.1in}
\caption{The number of iterations and CPU time of the standard, Jacobi, Gauss-Seidel,
and Gauss-Seidel-Jacobi value iteration algorithms with and without an accelerating
operator applied to a family of MDPs from Example~\ref{ex:diag}} \label{tab:2}
\end{table}
\end{center}

\end{document}